\documentclass[11pt,letterpaper]{article}
\usepackage[utf8]{inputenc}
\usepackage{amsmath}
\usepackage{amsfonts}
\usepackage{amssymb}
\usepackage{caption}
\usepackage{subcaption}
\usepackage{xcolor}
\usepackage{hyperref}
\usepackage{graphicx}
\usepackage{authblk}
\usepackage{cleveref}
\usepackage[letterpaper]{geometry}

\newcommand{\comment}[1]{}
\title{Exponential integrators for non-linear diffusion \\[1ex] \large LLNL-TR-832157}
\author[1]{Valentin Dallerit} 
\author[1]{Mayya Tokman}
\author[2]{Ilon Joseph}
\affil[1]{School of Natural Sciences, University of California, Merced, CA 95343, USA}
\affil[2]{Lawrence Livermore National Laboratory, P.O. Box 808, Livermore, CA 94551, U.S}
\begin{document}

\maketitle 

\tableofcontents
\pagebreak

\section{Introduction}

The goal of this project is to compare the performance of exponential time integrators with traditional methods such as diagonally implicit Runge-Kutta methods in the context of solving the system of reduced magnetohydrodynamics (RMHD). 
In this report, we present initial results of a proof of concept study that shows that exponential integrators can be an efficient alternative to traditional integration schemes. 

For this work the spatial discretization is done using the finite element method and we utilize LLNL's MFEM software package for this purpose. Discretizing RMHD equations in space using finite elements leads to a large system of ordinary differential equations (ODEs) of the form: 
\begin{equation}
\begin{cases}
\label{eq:ode}
y'(t) = f(y(t)) \\
y(t_0) = y_n
\end{cases} 
\end{equation} 
where $y(t)$ is the vector of degrees of freedom of the system and $t_0 \le t \le t_f$. 
Due to the presence of a widely ranging timescales in the RMHD equations, system (\ref{eq:ode}) is very stiff. Therefore, it is usually solved using an implicit time integrator, such as implicit Runge-Kutta or Backward Differentiation Formula (BDF) method. Each time step of an implicit method requires the solution of a large system of nonlinear equations. Approximating such a solution is computationally expensive, and a preconditioner is necessary. Often, preconditioners that yield sufficient efficiency are hard to construct. 
For these reasons, we want to explore the use of exponential integrators for this problem. These methods are attractive for large stiff problems since they have good stability properties, allow larger time step, and their cost per time step can be computationally cheaper than implicit methods.

In order to show the advantages of exponential time integrators, we formulate two nonlinear diffusion test problems corresponding to a simplified version of the RMHD model of interest. We then compare the efficiency of these exponential schemes to the traditional methods. 

This report is organized as follows. First, we introduce the numerical exponential integrators  used for this comparison. Then, we describe the numerical test problems and detail the implementation. Finally, we present the results of the comparison study and demonstrate that exponential time integrators can lead to a more accurate and more efficient solution compared to implicit schemes.

\section{Time integration} \label{sec:review}
\subsection{Time stepping methods}
In this section, we give a brief overview of exponential time integrators and present the schemes we will be using in the numerical experiments. We also provide details on the implementation and software we use.
 
Exponential integrators are usually derived in the following way. First, the right-hand side function $f$ of \cref{eq:ode} is expanded in a Taylor series around the solution $y(t_n)$ at a given time $t_n$:
\begin{equation}
  f(y) = f(y(t_n)) + J_n (y - y(t_n)) + R(y)
\end{equation}
where $J_n$ is the Jacobian matrix of $f$ at time $t_n$ and $R(y) = f(y) - f(y(t_n)) - J_n (y - y(t_n))$ is the nonlinear remainder of the Taylor expansion.
Using the integrating factor $e^{-J_n t}$ and integrating the equation from $t_n$ to $t_{n+1}$, we can get the solution at a future time $t_{n+1} = t_n + h$ as: 
\begin{align}
    \label{eq:beg}
    y(t_{n+1}) = y(t_n) + \varphi_1(hJ_n) h f(y(t_n)) + \int_0^1 e^{(1-\theta) h J_n} \ h R(y(t_n + h \theta)) \ d \theta
\end{align} 

where the exponential-like functions $\varphi_k$  are defined as:
\begin{align*}
\varphi_0(z) &= e^z = \sum_{n=0}^\infty \frac{z^n}{n!}\\ 
\varphi_k(z) &= \int_0^1 e^{(1-\theta) z }\frac{\theta^{k-1}}{(k-1)!}d\theta  = \sum_{n=0}^\infty \frac{z^n}{(n+k)!}  \ \ \text{for} \ \ k \ge 1
\end{align*}

\Cref{eq:beg} is the integral form of system (\ref{eq:ode}) and a starting point for the derivation of most exponential integrators. An exponential integrator is derived by choosing a numerical approximation of the nonlinear integral in \cref{eq:beg}.
For the study, we choose EPI type methods (\cite{tokmanEfficientIntegrationLarge2006}, \cite{tokmanNewClassExponential2011}) derived in \cite{rainwater2016new} as they have been shown to be efficient for this type of problems (\cite{tokmanNewAdaptiveExponential2012}, \cite{einkemmer2017a}).  
In particular, we will be testing the following two exponential schemes:

\begin{itemize}
\item EPI2 / Exponential Euler ($2^{nd}$ order):
  \begin{align}
    \label{eq:epi2}
    y_{n+1} &= y_n + \varphi_1(hJ_n) hf_n
  \end{align}
  \item EPIRK4 ($4^{th}$ order)  \cite{rainwater2016new} : 
  \begin{alignat}{3}\label{eq:epirk4}
  Y_{1} &= y_n  &+&  \frac{1}{8} \varphi_1\left(\frac{1}{8}h J_n\right)hf_n \notag \\
  Y_{2} &= y_n  &+&  \frac{1}{9}\varphi_1\left(\frac{1}{9} h J_n\right)hf_n \notag \\
  y_{n+1} &= y_n &+& \varphi_1\left( h J_n \right) h f_n \notag \\ 
  & & +& \varphi_3(h J_n) (\alpha_{3,1} hR(Y_1) + \alpha_{3,2} hR(Y_2)) \notag\\ 
  & & +& \varphi_4(h J_n) (\alpha_{4,1} hR(Y_1) + \alpha_{4,2} hR(Y_2))
\end{alignat}
\end{itemize}
where $\alpha_{3,1} = -1024, \ \alpha_{3,2} = 1458, \ \alpha_{4,1} =  27648, \ \alpha_{4,2} = -34992$.

Note that the  major computational cost of exponential integrator is the computation of the product between exponential functions of the jacobian matrix and vector. The EPI methods are specifically designed to reduce the number of such computations. For example, EPIRK4 requires only 2 such evaluations. As we explain bellow, these evaluations are usually done using Krylov projection type of technics so one would compare the number of such evaluations with the number of nonlinear iterations required at each time step of an implicit method.

For comparison we will use the following explicit and implicit methods implemented in MFEM:
\begin{itemize}
  \item Explicit methods: explicit Euler (\emph{ForwardEulerSolver}), Runge-Kutta of order 2, 3 and 4 (\emph{RK2Solver, RK3SSPSolver, RK4Solver})
  \item Implicit methods: implicit Euler (\emph{BackwardEulerSolver}), singly-diagonally implicit Runge-Kutta of order 2 and 3 (\emph{SDIRK23Solver, SDIRK23Solver})
\end{itemize}

\subsection{Exponential matrix function evaluation}

The main computational challenge in the time stepping of exponential schemes is the evaluation of the matrix-vector products involving the $\varphi_k$ functions. 
For the problems we are considering, the Jacobian matrix involved in these computations is large and stiff. Therefore, we need an efficient method adapted to this kind of problems. In this work, we are using the KIOPS algorithm \cite{gaudreaultKIOPSFastAdaptive2018}. It is using a Krylov subspace to project the large operator into a smaller space where the computation can be carried more efficiently. This idea is similar to the one used in methods like GMRES or conjugate gradient.
Moreover, some optimizations are used in this method to improve the efficiency. 
First, using a theorem from \cite{al-mohyComputingActionMatrix2011a}, it is possible to reduce the computation of a linear combination of $\varphi$ functions of the form $\sum_{i=0}^k \varphi_i(hJ) v_i$ to the computation of a single matrix exponential times a vector.
Then, in \cite{gaudreaultKIOPSFastAdaptive2018}, the authors show that during the Arnoldi process it is enough to orthogonalize the new vector with respect to the last 2 previous vectors. This modification reduces the number of dot products from $O(m^2)$ to $O(m)$, where $m$ is the size of the Krylov space, while keeping the procedure stable. This is especially important for parallel code, as dot products require a communication between the computation nodes.
Unlike standard Arnoldi procedure and adaptive Krylov algorithm, KIOPS requires only 2 dot products per Krylov vector which significantly improve it's parallel efficiency.
The last optimization, originally presented in \cite{niesen2012a}, is to substep the computation
of the matrix exponential by using the following property of the exponential function:
$$ e^A  v = e^{\tau_k A} ... e^{\tau_2 A}e^{\tau_1 A} v \ \ \text{if} \ \ \sum_{i=1}^k \tau_i = 1 $$
Using this equation, it is possible to compute $e^A v$ by first computing $w_1 = e^{\tau_1 A} v$, then computing $w_2 =
e^{\tau_2 A} w_1$ , and so on until $\tau_k$. By choosing $\tau_i < 1$, we have $||\tau_i A|| < ||A||$ and therefore the approximation requires a smaller Krylov subspace. This also allows the computation of $e^{\tau A}$ with several values of $\tau$ with a single Krylov projection. KIOPS uses an algorithm to compute the values of $\tau$ adaptively.

Using these optimizations, it is possible to advance the solution using the method EPIRK4 with only 2 Krylov projections instead of the 5 projections require with a naive implementation. To do so, we first compute the stages $Y_1$ and $Y_2$ in a single computation using the substepping. Then, we can compute the linear combination of $\varphi$ functions with a second projection using the first optimization.
For full details on EPI methods with KIOPS see \cite{gaudreaultKIOPSFastAdaptive2018} and \cite{gaudreault2022high}. 

\section{Implementation}
The KIOPS method as well as the exponential integrators are implemented in the \emph{EPIC-cpp} package available at \url{https://gitlab.com/tokman-lab/epic-cpp}. This package only requires the user to provide the right-hand side function and, optionally, the associated jacobian matrix or a matrix-free implementation of it. It can be used for both serial and parallel applications (using the MPI standard) and is built on top of \emph{NVector} module of \emph{SUNDIALS} for vector storage and basic operations.

As part of this project, an interface between \emph{EPIC-cpp} and \emph{MFEM} was developped. This allows users to access and use the integrators available in the \emph{EPIC-cpp} package directly from \emph{MFEM}. This wrapper works similarly to the integration of SUNDIALS into MFEM and is transparent for the user. One can use the \emph{EPIRK4} or \emph{EPI2} class in place of the other MFEM ode solver (e.g. \emph{ForwardEulerSolver}, \emph{SDIRK23Solver}) and call the same methods. Additional parameters can be used during the construction of the exponential integrator to choose between user-provided jacobian or finite difference approximation. 
This interface is currently available in the \emph{epic-dev} branch of MFEM. 
The other time integration methods used for the comparison are from the already available MFEM ode solver.
The test problems used in the numerical experiments are implemented in the \emph{feature/exponential-integrator} branch of the \emph{mhdex} project, available at \url{https://lc.llnl.gov/bitbucket/scm/mwm/mhdex.git}.
All the necessary dependancies and test problems can be obtained from the \emph{feature/exponential-integrator} branch of the \emph{mdh-mfem} project, available at \url{https://lc.llnl.gov/bitbucket/scm/mwm/mhd-mfem.git}.

\section{Test problems} \label{sec:experiment}
In order to compare the performance of the different time integration methods, we set up two test problems. 
These problems model nonlinear diffusion in a way that is similar to what is found in the RMHD problem of interest.

\begin{itemize}
    \item 1D diffusion problem:  

      The first problem correspond to a nonlinear diffusion PDE in 1 dimension. 
    \begin{align}
      \label{eq:1d_og}
      \frac{\partial u}{\partial t} = \frac{\partial}{\partial x} \left( \mu(u) \frac{\partial u}{\partial x} \right) + s(x)
    \end{align}
on the domain $x \in [0,1]$
with $\mu(u) = (\beta_1 + \beta_2 u^{5/2})$ and the source is defined by $s(x) = e^{ -\frac{1}{2} ( (x-1/2)/\sigma )^2 }$.
The diffusion coefficients can be selected in order to adapt the stiffness and the amount of nonlinearity of the problem. With $\beta_2 = 0$, the problem is linear. In the other case, the nonlinearity is chosen to emulate the kind of diffusivity found in RMHD problems. 

We found that in order to get the correct order of convergence, we need to rewrite the problem in the following form:
    \begin{align} 
      \label{eq:1d_mod}
      \frac{\partial u}{\partial t} = \frac{\partial^2}{\partial x^2} \left( g(u) \right) + s(x)
    \end{align}

    where $g(u) = \beta_1 u + \frac{2}{7} \beta_2 u^{7/2} = \int \mu(u) du$. Since mathematically \cref{eq:1d_og} and \cref{eq:1d_mod} are equivalent, the discrepency is most likely due to a bug in the current version  of MFEM.

\item 2D anisotropic diffusion

  The second problem is an anisotropic diffusion PDE in 2 dimensions. $$ \frac{\partial u}{\partial t} = \nabla . \left( \mu(u) \nabla u \right) + s(x)$$

on the domain $(x,y) \in [0,1]^2$
    with $\mu(u) = \kappa  \left[(\beta_1 + \beta_2 u^{5/2}) (\hat{b} \otimes \hat{b}) + 10^\alpha (I - \hat{b} \otimes \hat{b})\right]$ and the source is defined by $s(x) = e^{ -\frac{1}{2} ( (x-1/2)/\sigma )^2 }$. The magnetic field $b$ used in this experiment is the \emph{2-wire} model and is represented in Figure \ref{fig:vectorfield} . 

\begin{figure}[!ht]
  \centering
  \includegraphics[scale=.4]{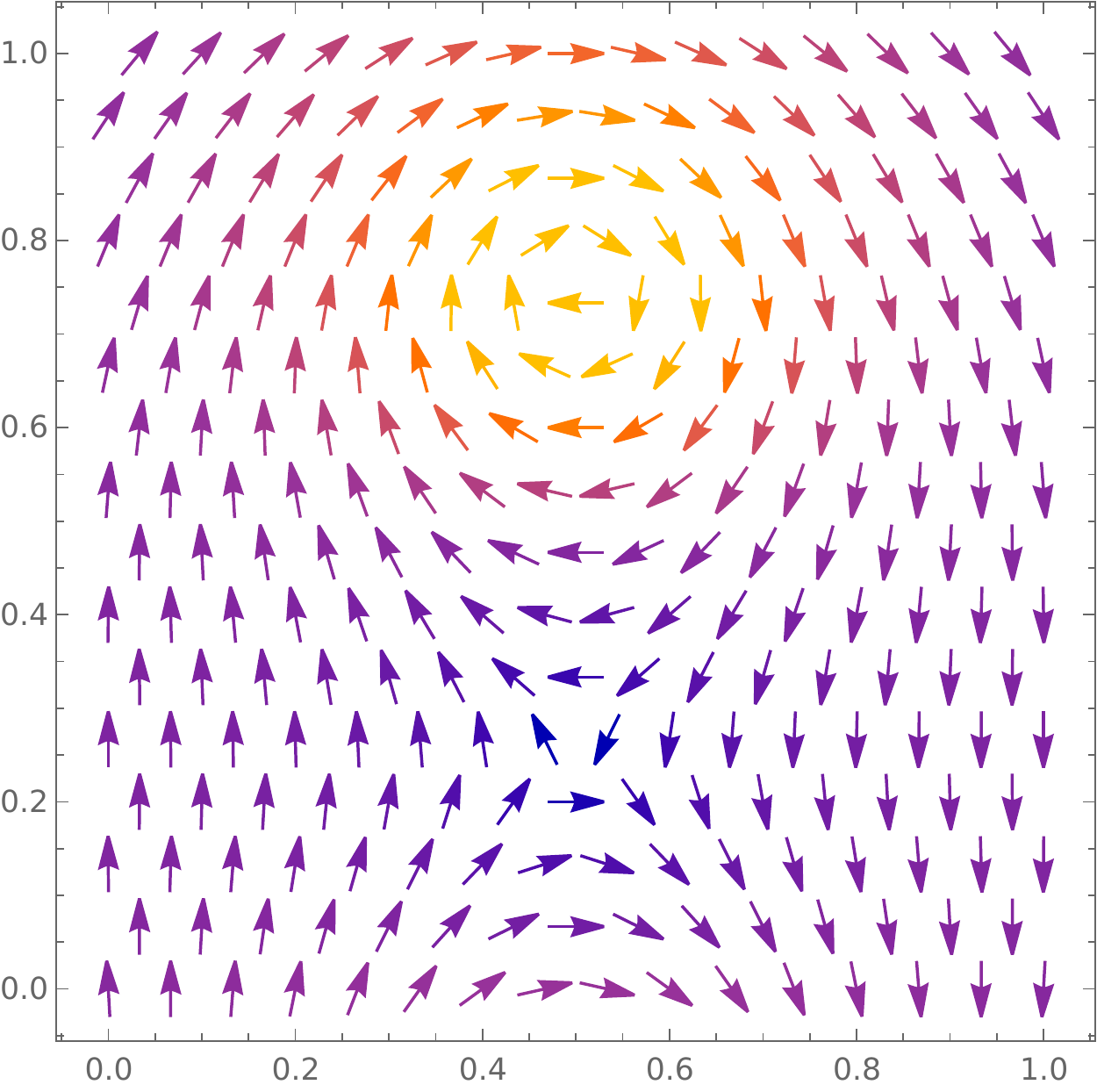}
  \caption{Magnetic field used for the anisotropic diffusion}
  \label{fig:vectorfield}
\end{figure}

    The $\beta$ parameters are controlling the diffusion in the direction of the magnetic field. As in problem (\ref{eq:1d_og}), it is possible to choose either a linear or nonlinear diffusion. The $\alpha$ coefficient is controlling the strength of the diffusion in the direction perpendicular to the magnetic field. As this parameter goes to $0$, a boundary layer is forming, making the problem stiffer. 

Similar to the 1D diffusion test, we need to rewrite the problem in the following form in order to get the correct order of convergence:
$$ \frac{\partial u}{\partial t} = \nabla . \left( \left(\hat{b} \otimes \hat{b} \right) \nabla \left(\beta_1 u + \frac{2}{7} \beta_2 u^{7/2} \right) \right) + \nabla . \left( (I - \hat{b} \otimes \hat{b}) \nabla \left(10^a u \right) \right) + s(x)$$

\end{itemize}

\clearpage

\section{Results} \label{sec:results}

\subsection{Validation of the order of convergence}
\begin{figure}[!ht]
  \centering
  \begin{subfigure}{.45\linewidth}
    \includegraphics[width=\linewidth]{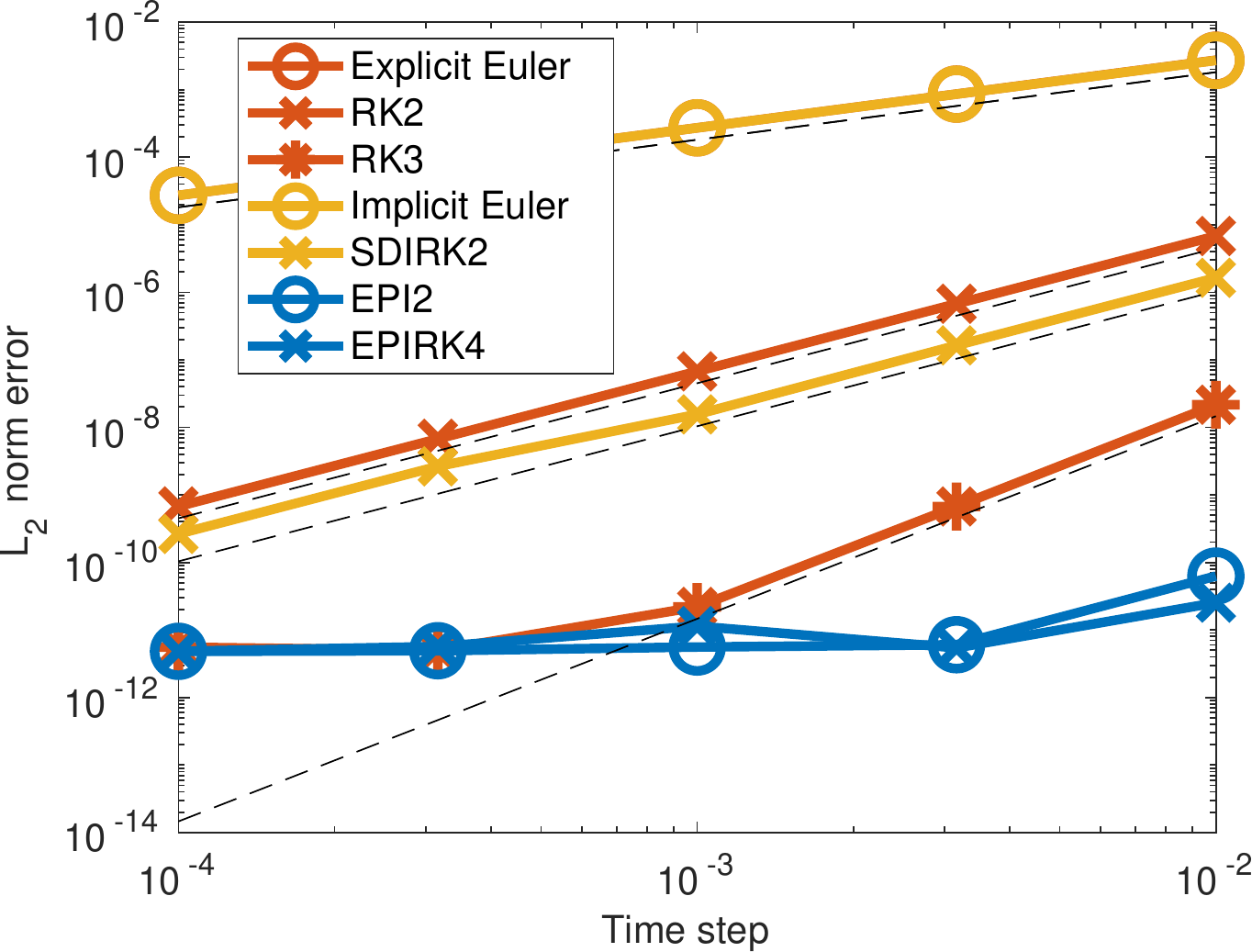}
    \caption{1d test problem \\ ($\beta_1 = 5 \times 10^{-3}, \beta_2 = 0$)}
  \end{subfigure}
  \quad
  \begin{subfigure}{.45\linewidth}
    \includegraphics[width=\linewidth]{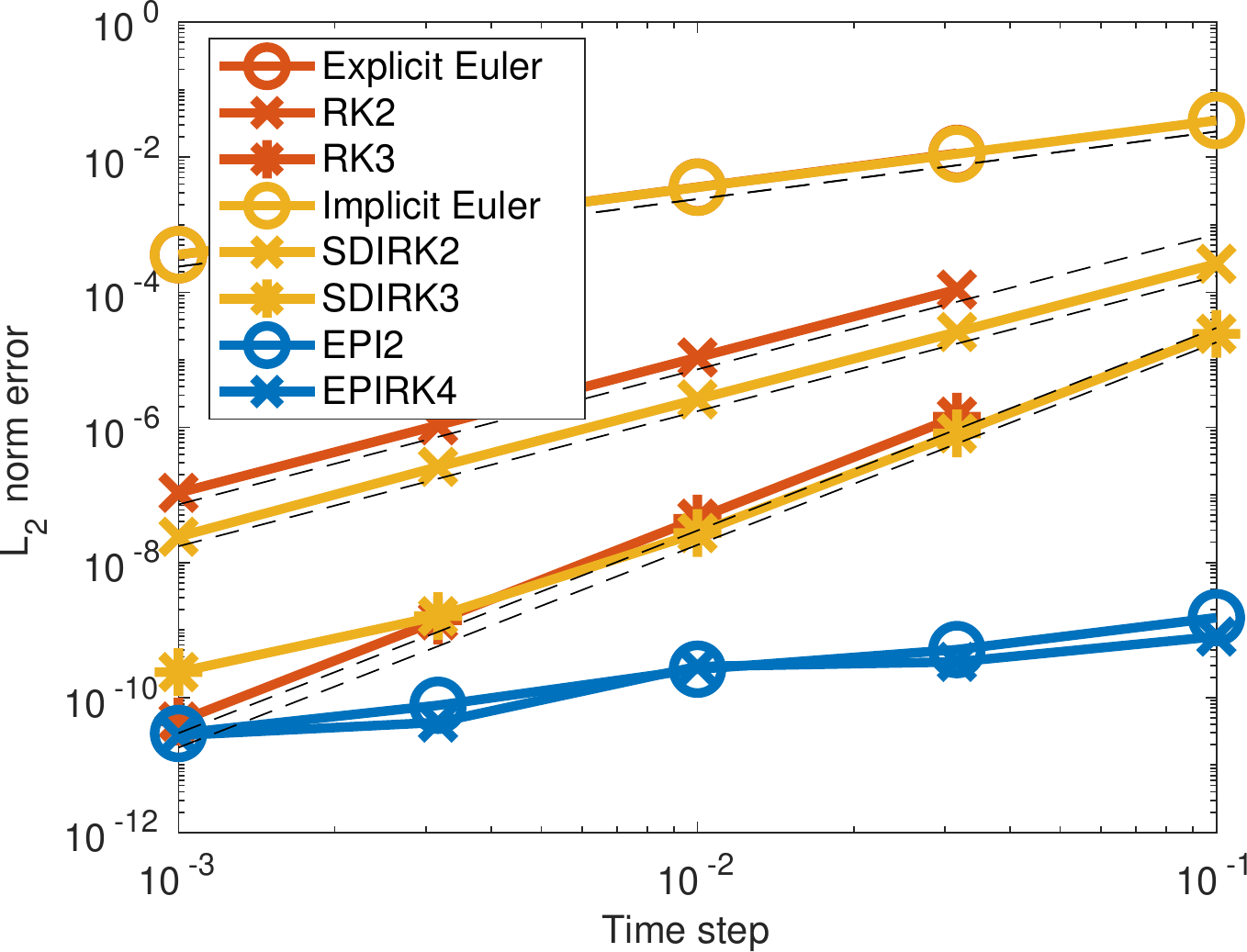}
    \caption{2d test problem  \\ ($\kappa = 10^{-2}, \alpha = 10^{-3}, \beta_1 = 1, \beta_2=0$)}
  \end{subfigure}
  \caption{Convergence diagrams for the 1d and 2d test problems with linear diffusion.}
  \label{fig:conv_lin}
\end{figure}

Our first task is to validate the corect order of convergence of the implemented methods. 
Since for linear problems, the schemes (\ref{eq:epi2}) and (\ref{eq:epirk4}) are exact (assuming that the $\varphi_k$ functions are evaluated exactly), we want to verify that if we set the nonlinear diffusion coefficient to zero, we are getting the expected result. \Cref{fig:conv_lin} shows the norm of the error as a function of the timestep for the 1d test problem on the left and the 2d test problem on the right with linear diffusion. As expected, for the non-exponential methods (red and yellow curves), the error is growing with the time step. For the exponential schemes, for all values of the time step, the error stays very low. This error is nonzero because of the tolerance in the evaluation of the $\varphi_k$ functions.

\begin{figure}[!ht]
  \centering
  \begin{subfigure}{.45\linewidth}
    \includegraphics[width=\linewidth]{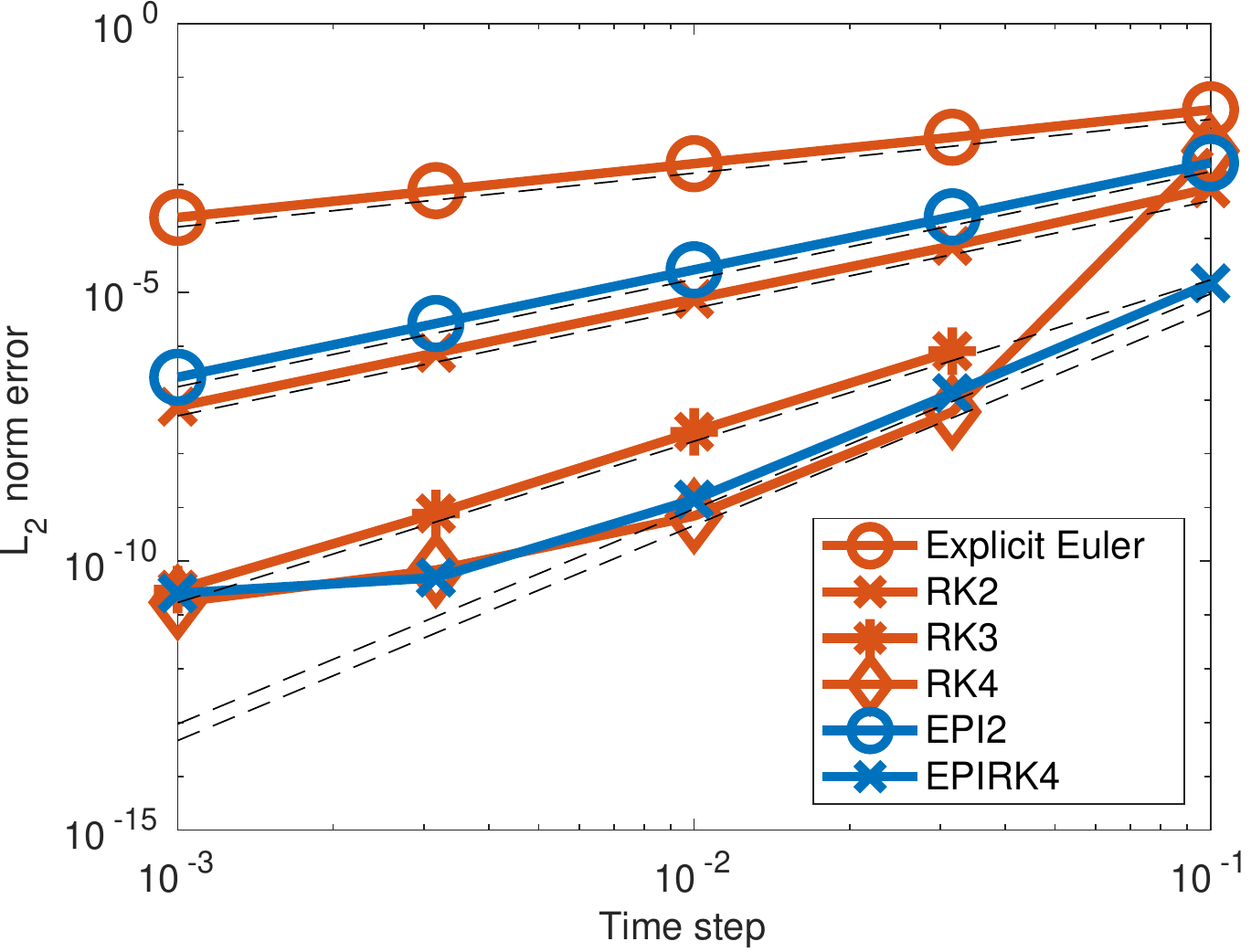}
    \caption{Explicit methods - Coarse grid ($n_{elem} = 50$)}
  \end{subfigure}
  \quad
  \begin{subfigure}{.45\linewidth}
    \includegraphics[width=\linewidth]{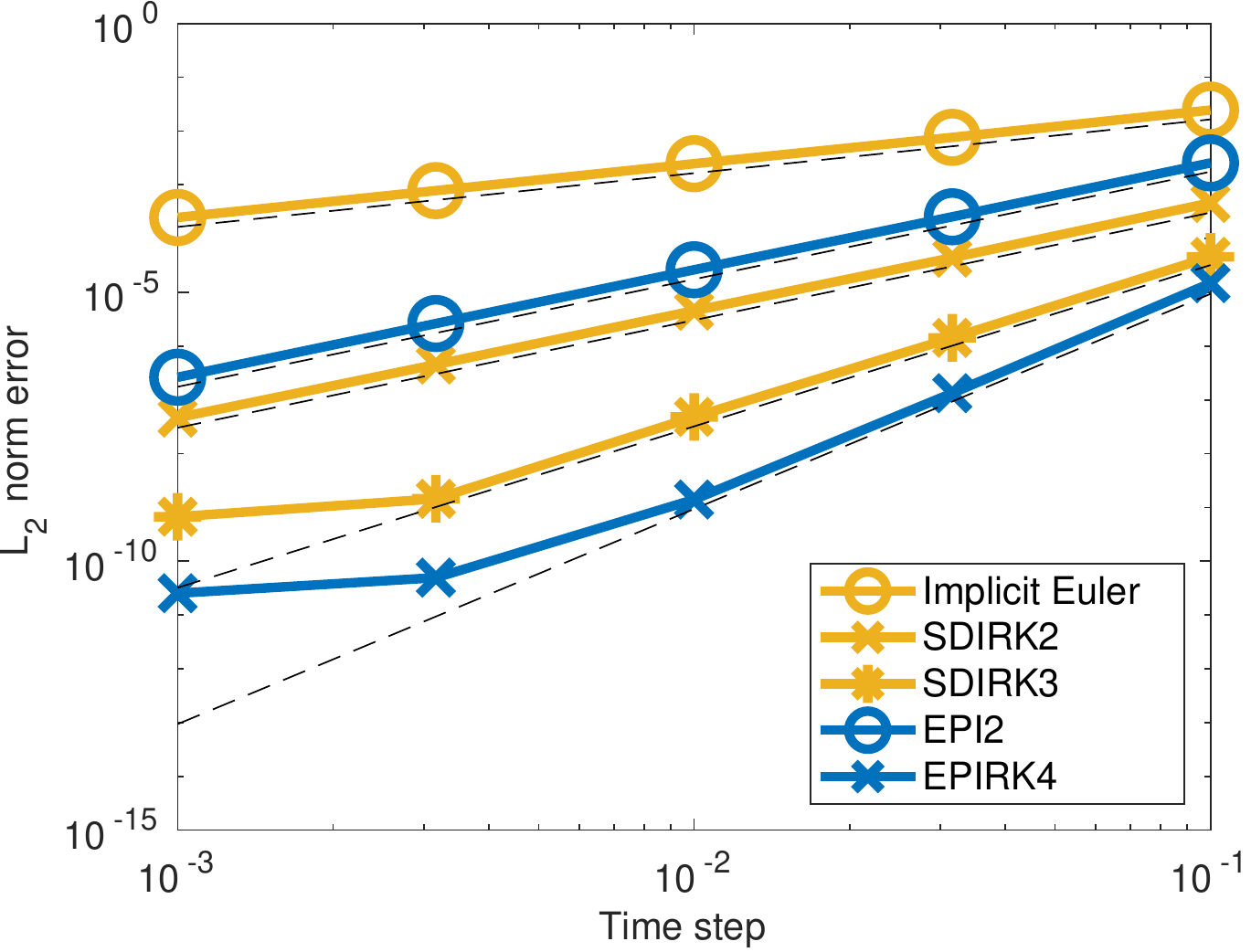}
    \caption{Implicit methods - Coarse grid ($n_{elem} = 50$)}
  \end{subfigure}
  \begin{subfigure}{.45\linewidth}
    \includegraphics[width=\linewidth]{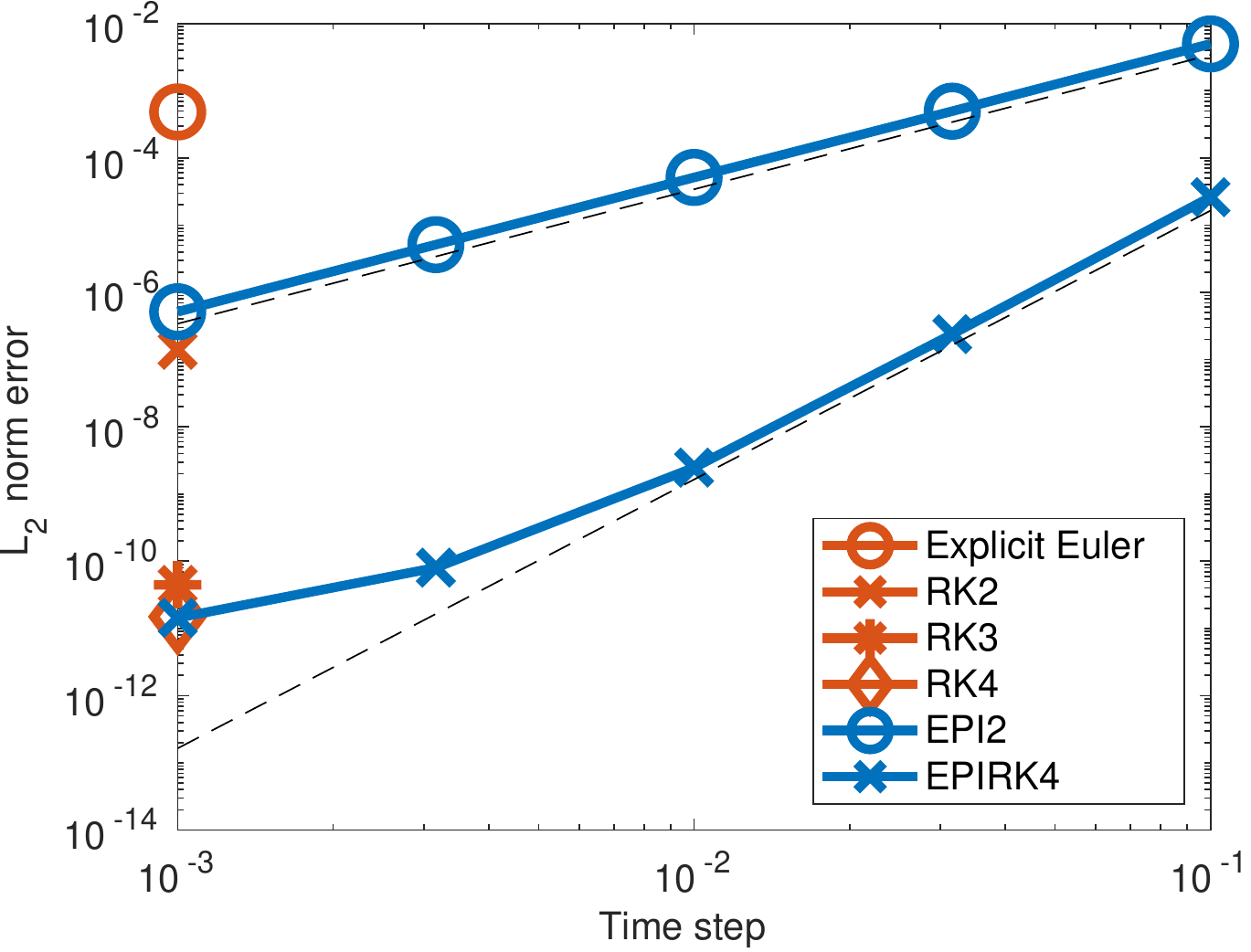}
    \caption{Explicit methods - Fine grid ($n_{elem} = 200$)}
  \end{subfigure}
  \quad
  \begin{subfigure}{.45\linewidth}
    \includegraphics[width=\linewidth]{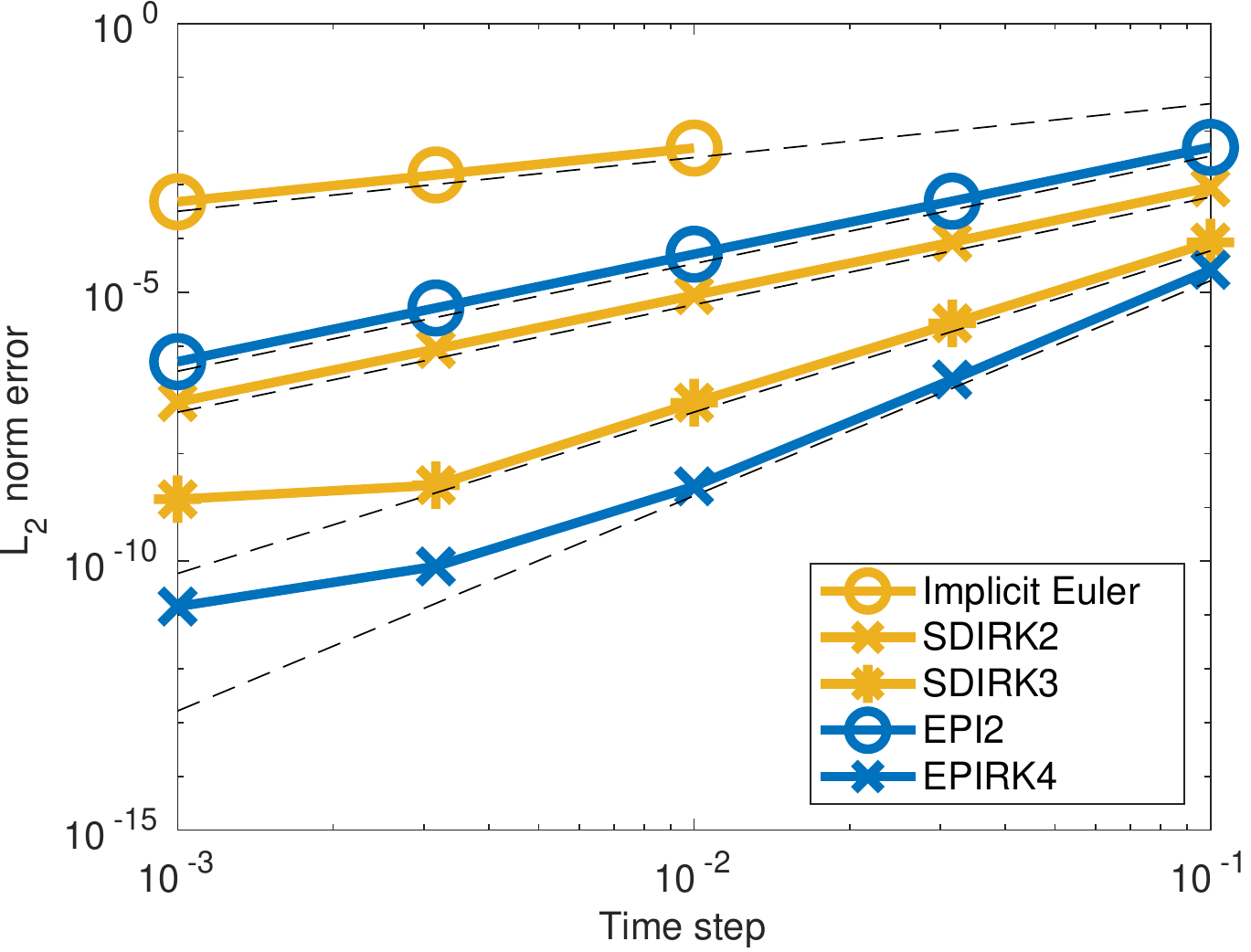}
    \caption{Implicit methods -  Fine grid ($n_{elem} = 200$)}
  \end{subfigure}

  \caption{Convergence diagrams for the 1d test problems with nonlinear diffusion ($\beta_1 = 5 \times 10^{-5}, \beta_2 = 5 \times 10^{-3}$).}
  \label{fig:conv_1d}
\end{figure}

\begin{figure}[!ht]
  \centering
  \begin{subfigure}{.45\linewidth}
    \includegraphics[width=\linewidth]{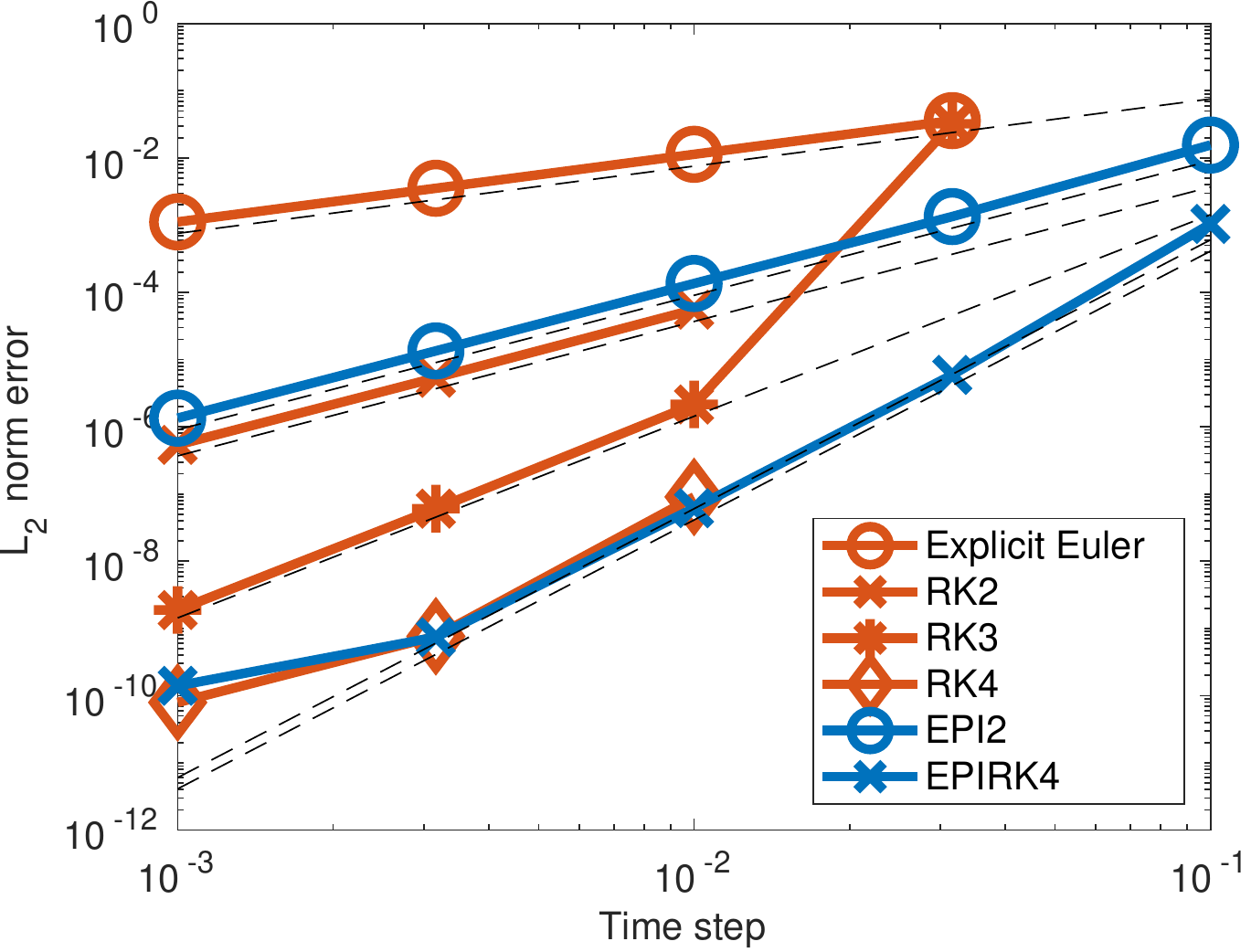}
    \caption{Explicit methods - Coarse grid ($n_{elem} = 20^2$)}
  \end{subfigure}
  \quad
  \begin{subfigure}{.45\linewidth}
    \includegraphics[width=\linewidth]{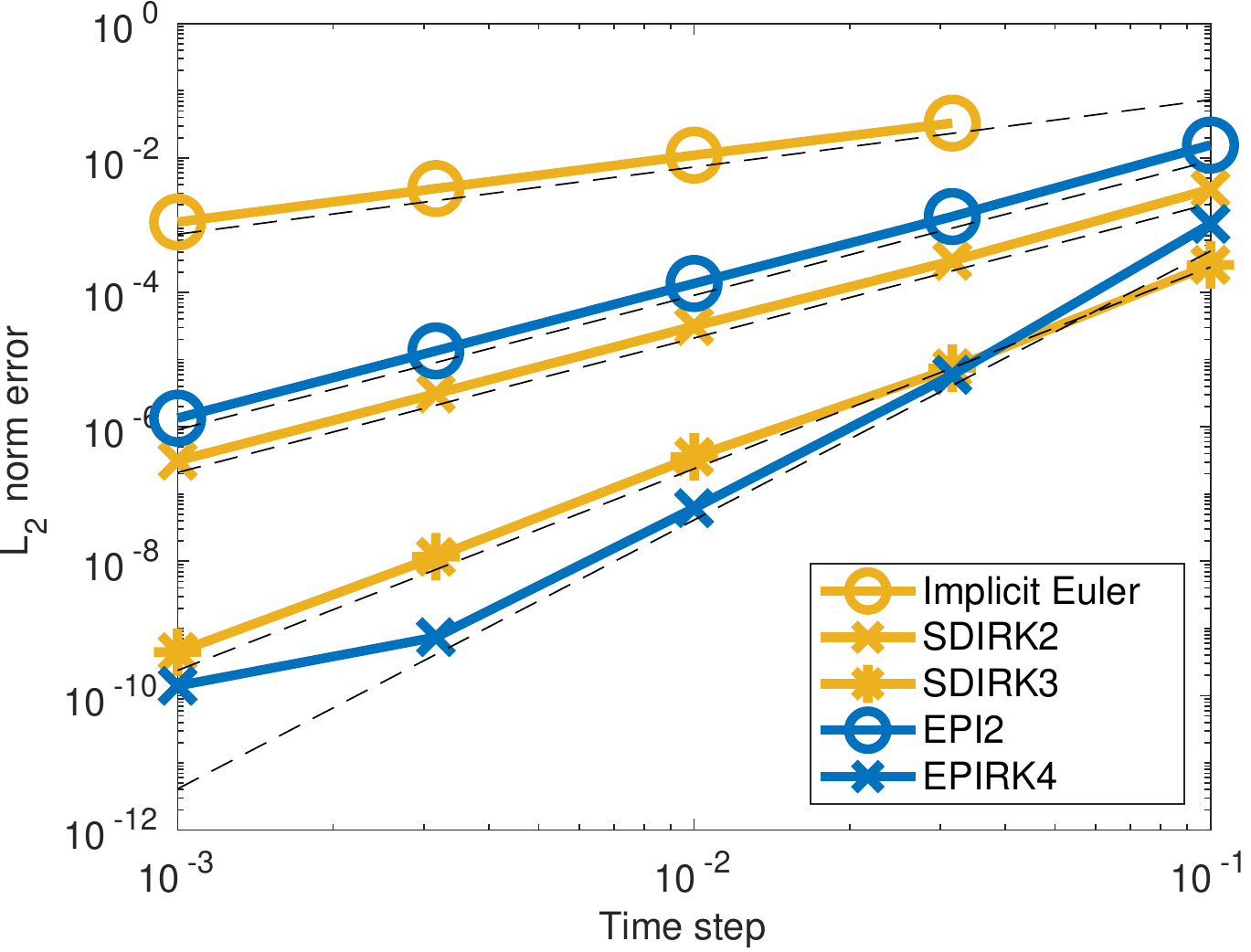}
    \caption{Implicit methods - Coarse grid ($n_{elem} = 20^2$)}
  \end{subfigure}
  \begin{subfigure}{.45\linewidth}
    \includegraphics[width=\linewidth]{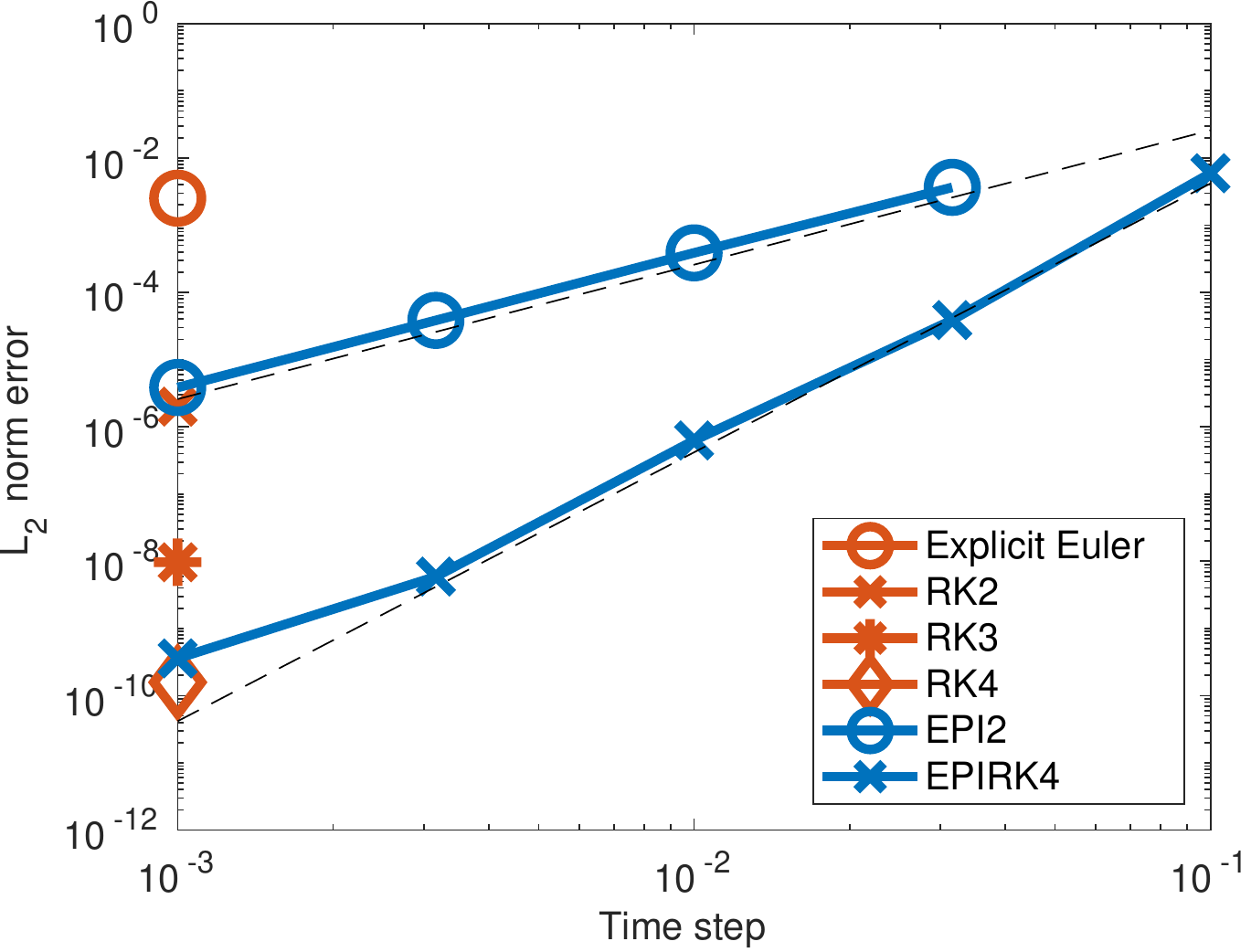}
    \caption{Explicit methods - Fine grid ($n_{elem} = 50^2$)}
  \end{subfigure}
  \quad
  \begin{subfigure}{.45\linewidth}
    \includegraphics[width=\linewidth]{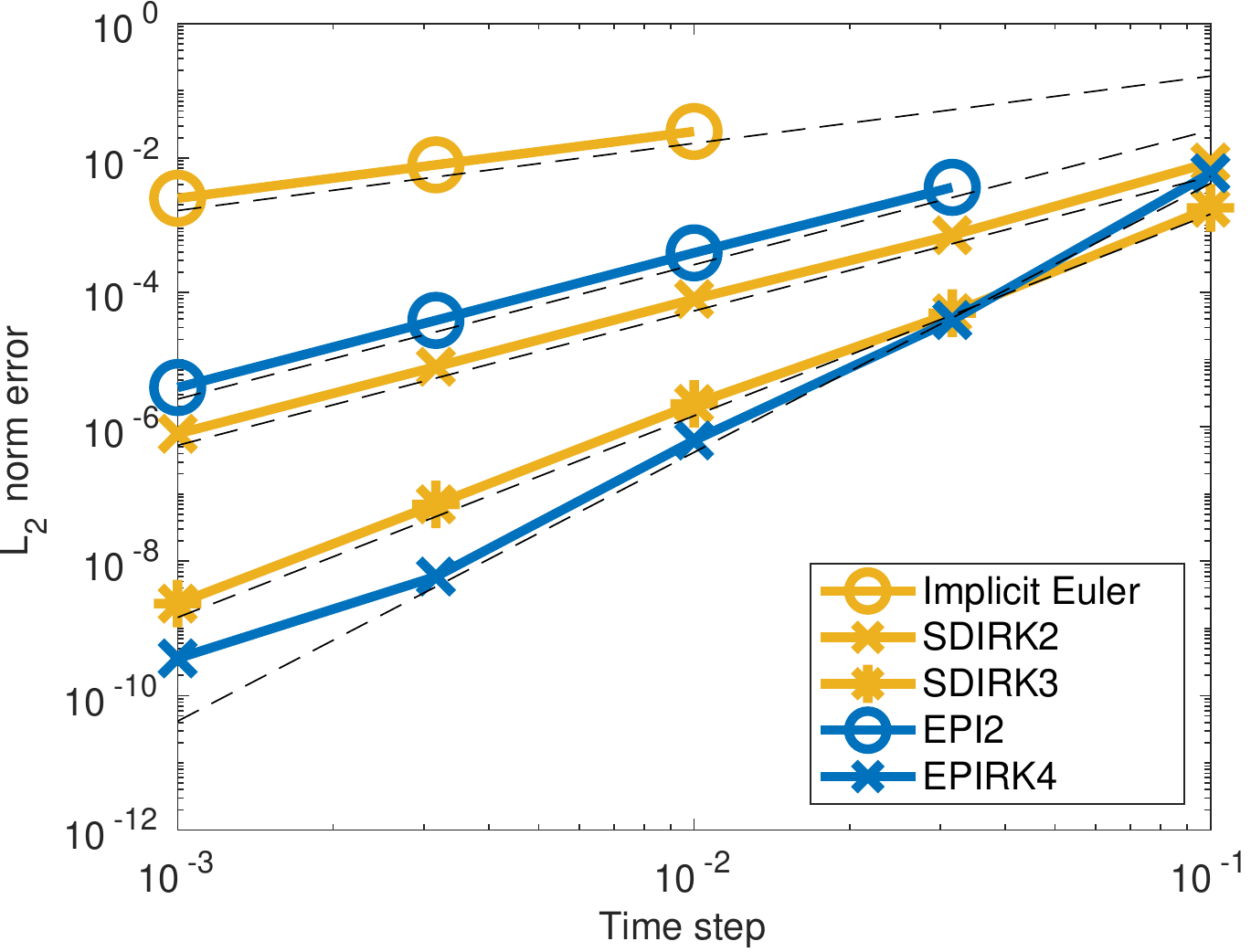}
    \caption{Implicit methods - Fine grid ($n_{elem} = 50^2$)}
  \end{subfigure}

  \caption{Convergence diagrams for the 2d test problems with anisotropic diffusion ($\kappa = 10^{-2}, \alpha = 10^{-3}, \beta_1 = 0, \beta_2=10$).}
  \label{fig:conv_2d}
\end{figure}

\Cref{fig:conv_1d} shows the convergence plot for the 1d test problem with nonlinear diffusion. The plots on the top correspond to a coarse grid making the problem non-stiff while the plots on the bottom are on a more refined grid making the problem stiffer. The left plots are comparing the explicit and exponential methods while the plots on the right are comparing the implicit and exponential methods. The dashed lines next to each curve correspond to the expected order of convergence. We can see that all the methods are converging at the expected rate. However, the explicit methods on the fine grid are only stable for the smallest time step while both the implicit and exponential methods do not have any stability issues. \Cref{fig:conv_2d} present the same results but for the 2d anisotropic diffusion test problem. With this problem, all methods still converge as expected and the explicit schemes still suffer from unstability on the finer grid.

\clearpage

\subsection{Performance comparison}
\begin{figure}[!ht]
  \centering
  \begin{subfigure}{.45\linewidth}
    \includegraphics[width=\linewidth]{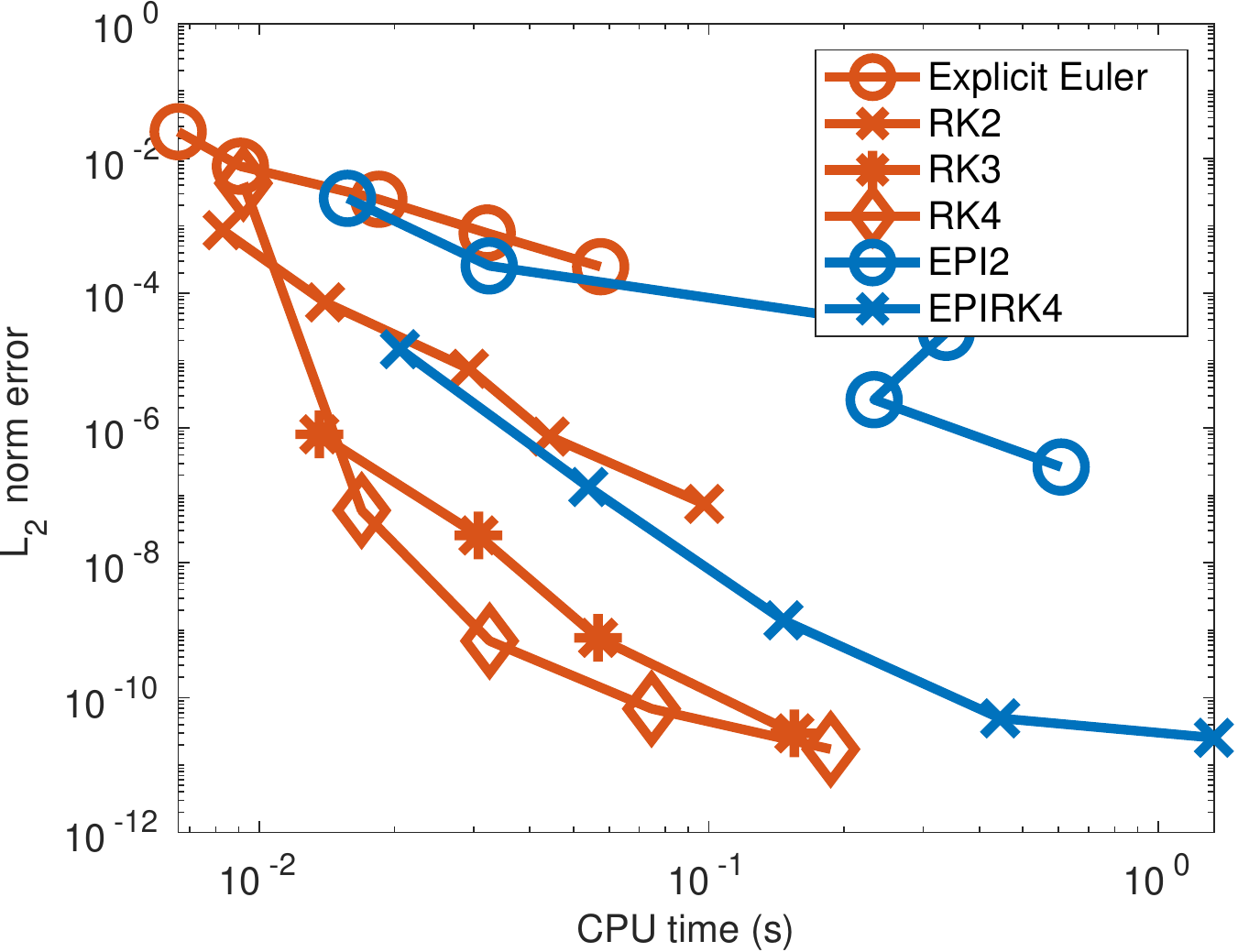}
    \caption{Explicit methods - Coarse grid ($n_{elem} = 50$)}
  \end{subfigure}
  \quad
  \begin{subfigure}{.45\linewidth}
    \includegraphics[width=\linewidth]{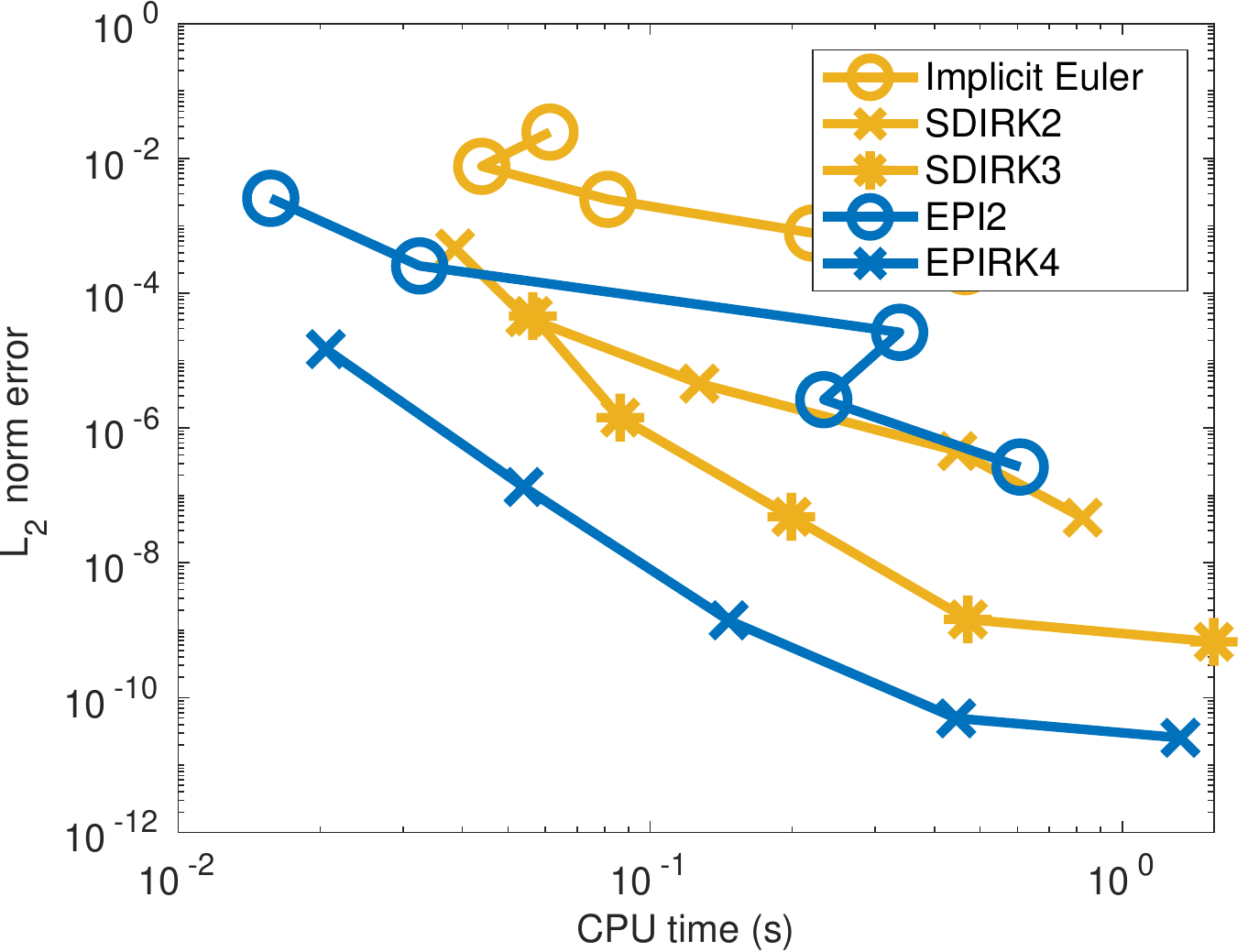}
    \caption{Implicit methods - Coarse grid ($n_{elem} = 50$)}
  \end{subfigure}
  \begin{subfigure}{.45\linewidth}
    \includegraphics[width=\linewidth]{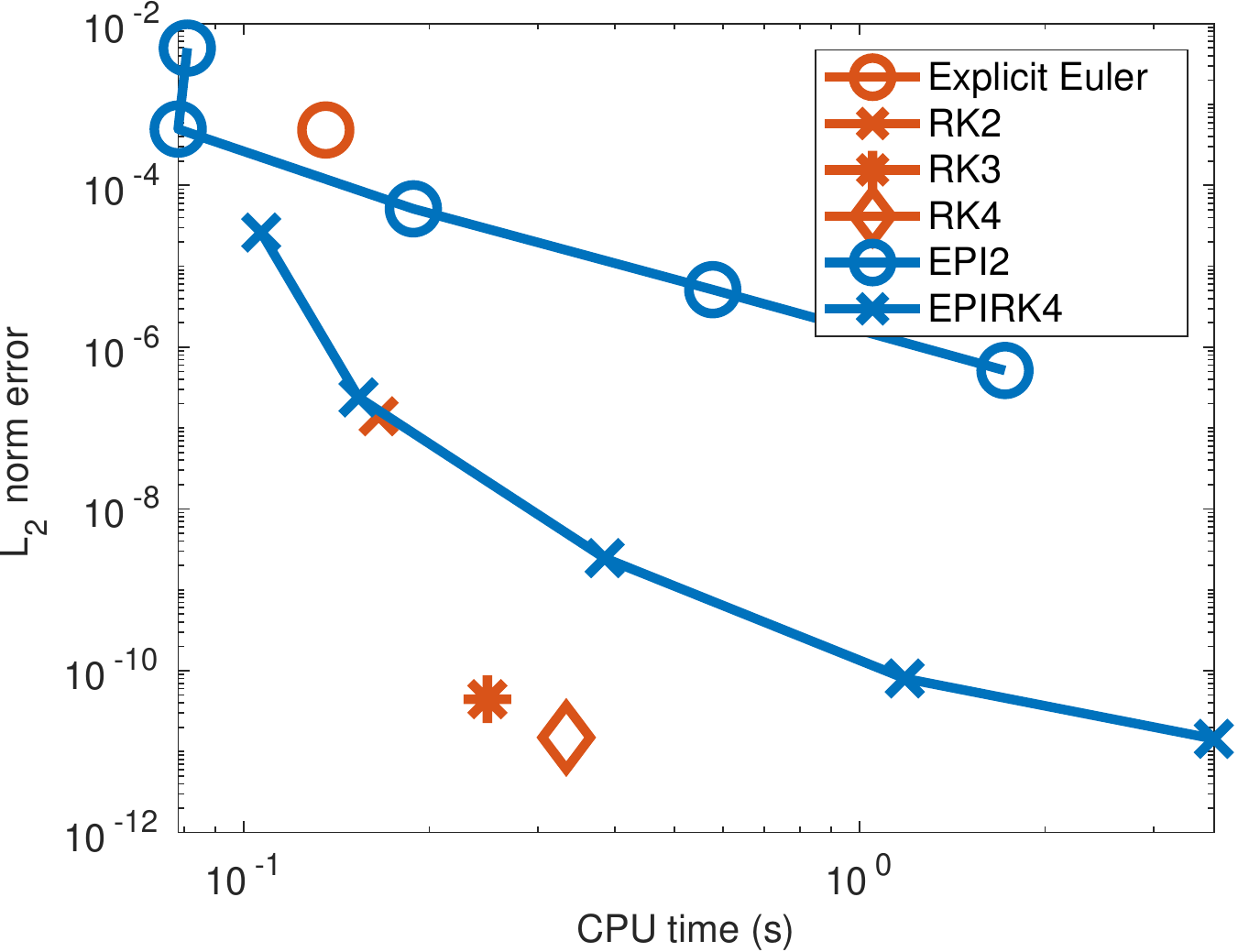}
    \caption{Explicit methods - Fine grid ($n_{elem} = 200$)}
  \end{subfigure}
  \quad
  \begin{subfigure}{.45\linewidth}
    \includegraphics[width=\linewidth]{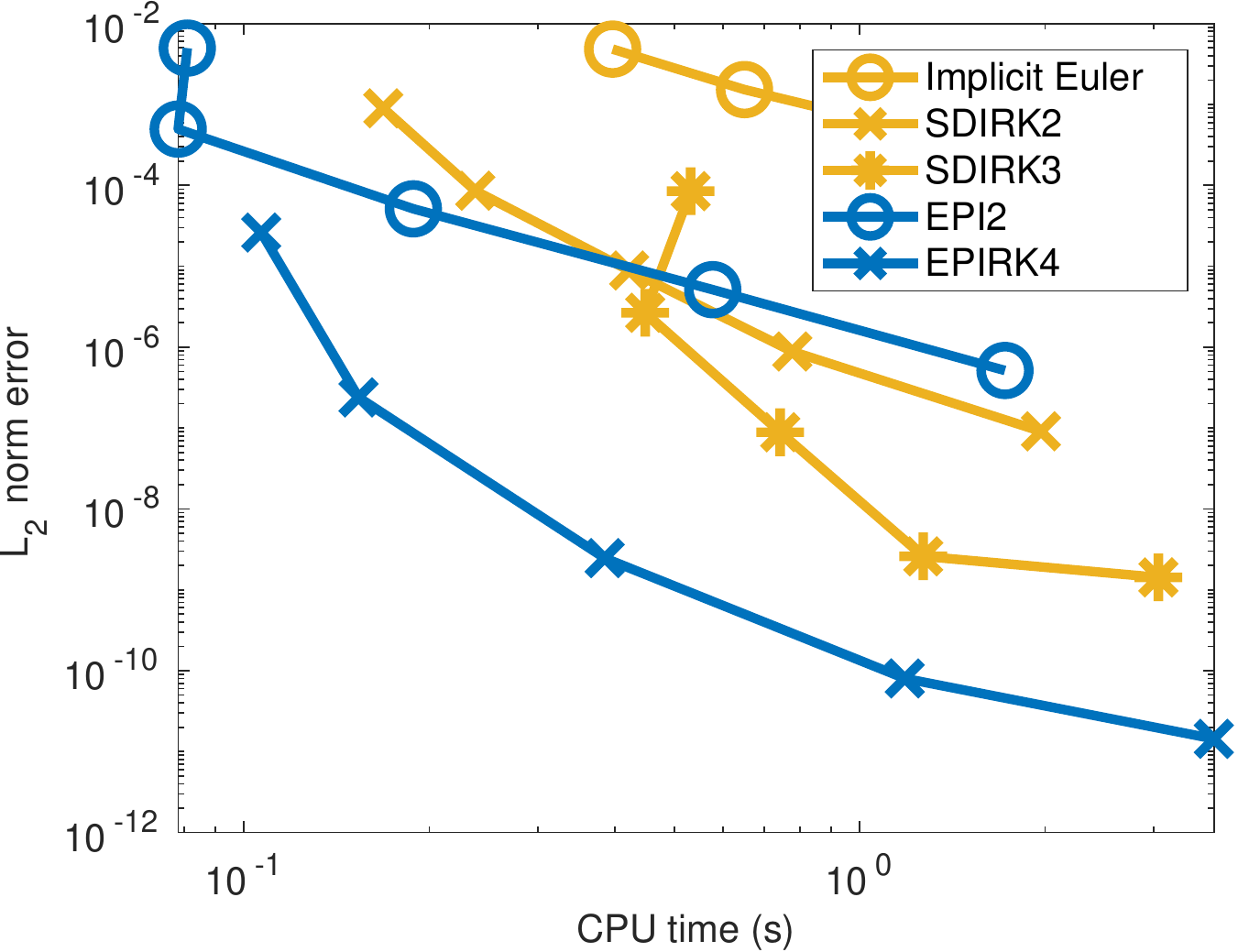}
    \caption{Implicit methods -  Fine grid ($n_{elem} = 200$)}
  \end{subfigure}

  \caption{Precision diagrams for the 1d test problems with nonlinear diffusion ($\beta_1 = 5 \times 10^{-5}, \beta_2 = 5 \times 10^{-3}$).}
  \label{fig:prec_1d}
\end{figure}

\begin{figure}[!ht]
  \centering
  \begin{subfigure}{.45\linewidth}
    \includegraphics[width=\linewidth]{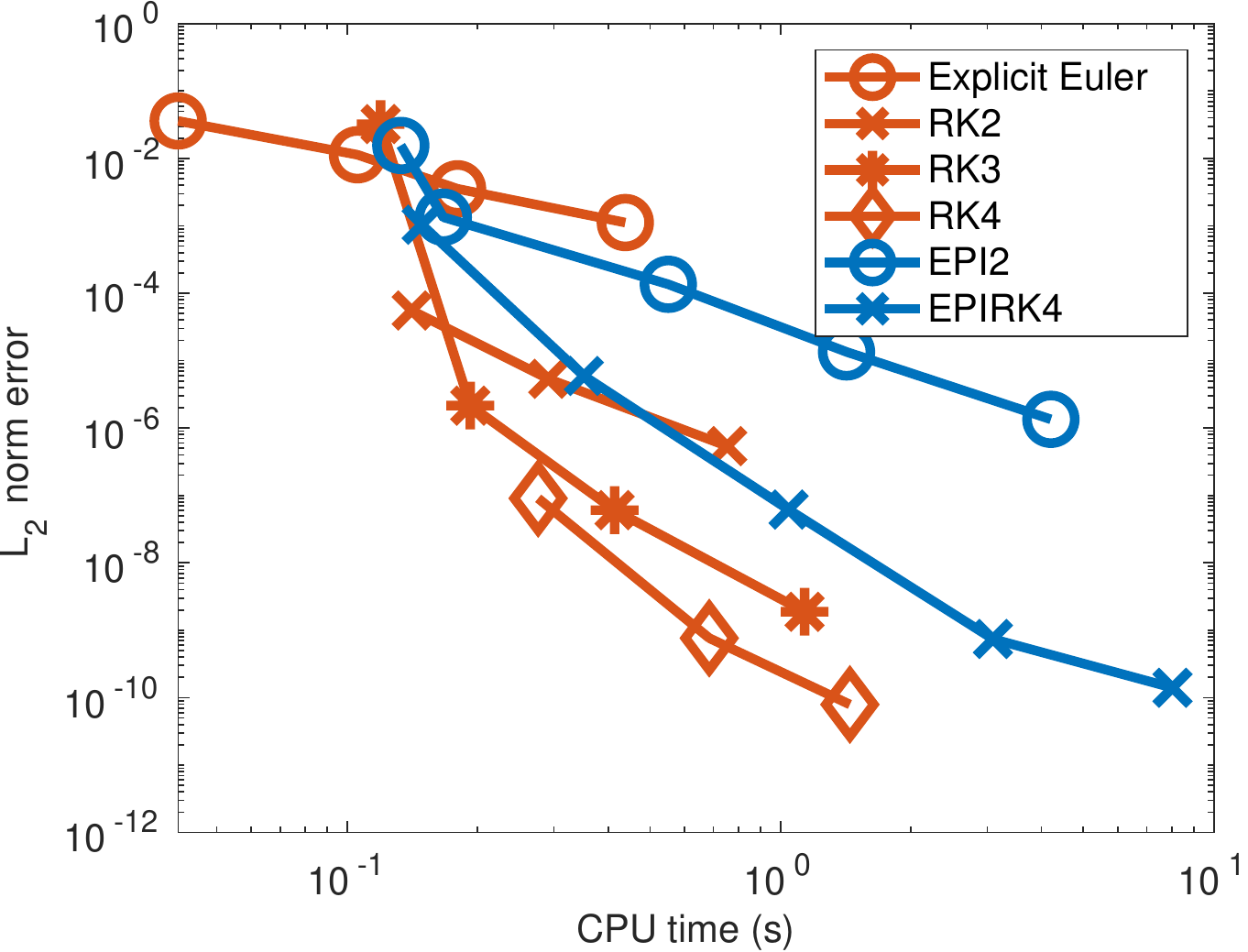}
    \caption{Explicit methods - Coarse grid ($n_{elem} = 20^2$)}
  \end{subfigure}
  \quad
  \begin{subfigure}{.45\linewidth}
    \includegraphics[width=\linewidth]{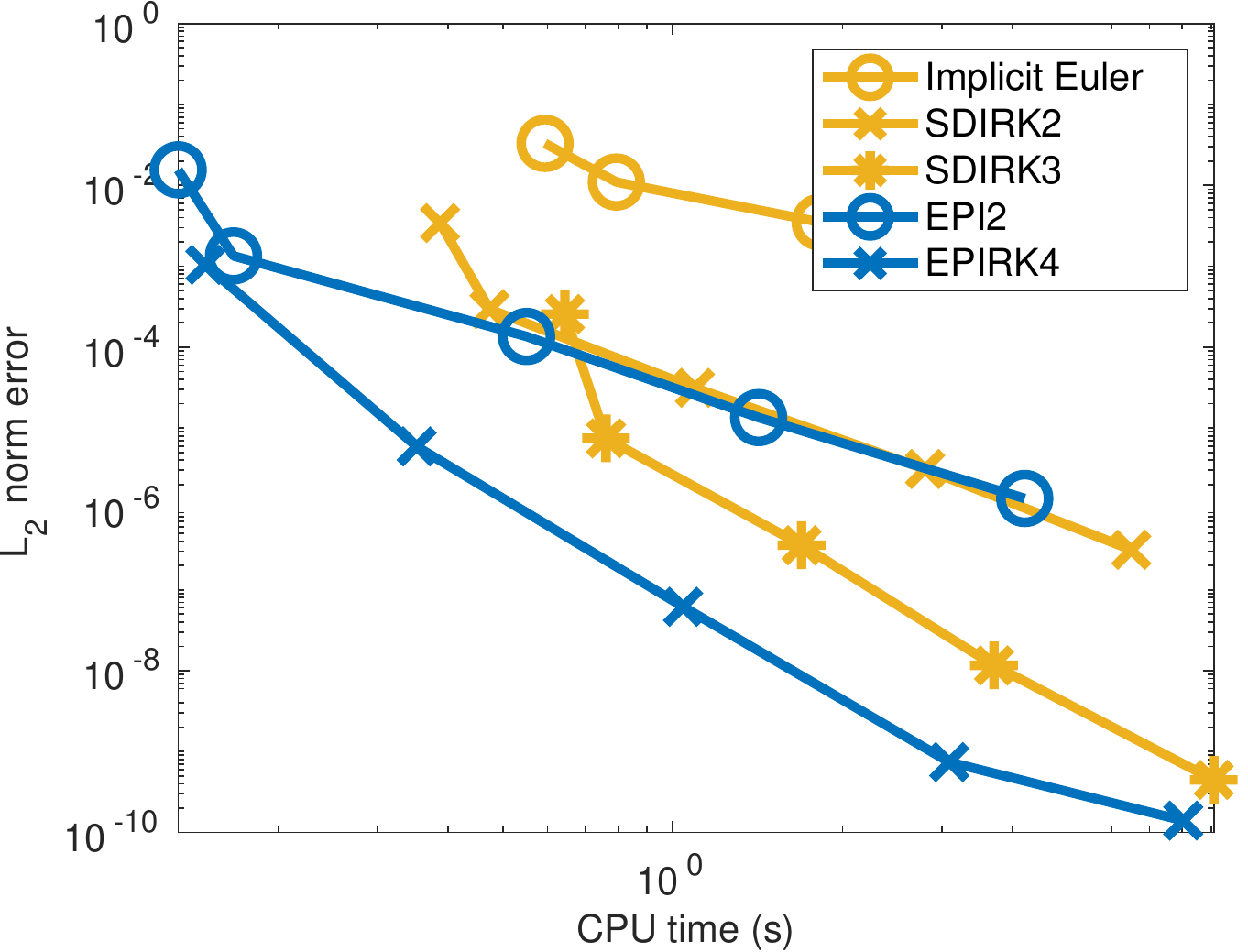}
    \caption{Implicit methods - Coarse grid ($n_{elem} = 20^2$)}
  \end{subfigure}
  \begin{subfigure}{.45\linewidth}
    \includegraphics[width=\linewidth]{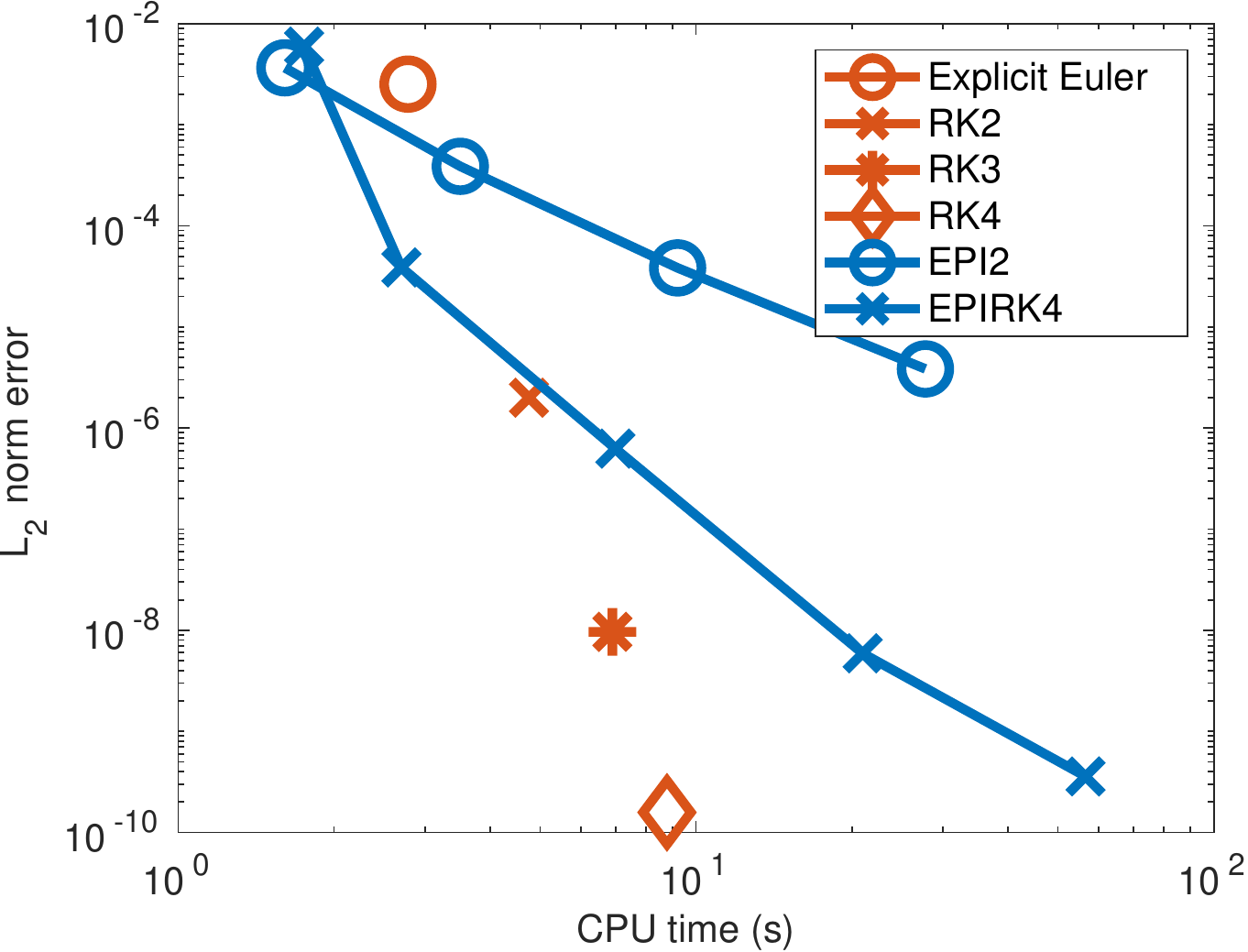}
    \caption{Explicit methods - Fine grid ($n_{elem} = 50^2$)}
  \end{subfigure}
  \quad
  \begin{subfigure}{.45\linewidth}
    \includegraphics[width=\linewidth]{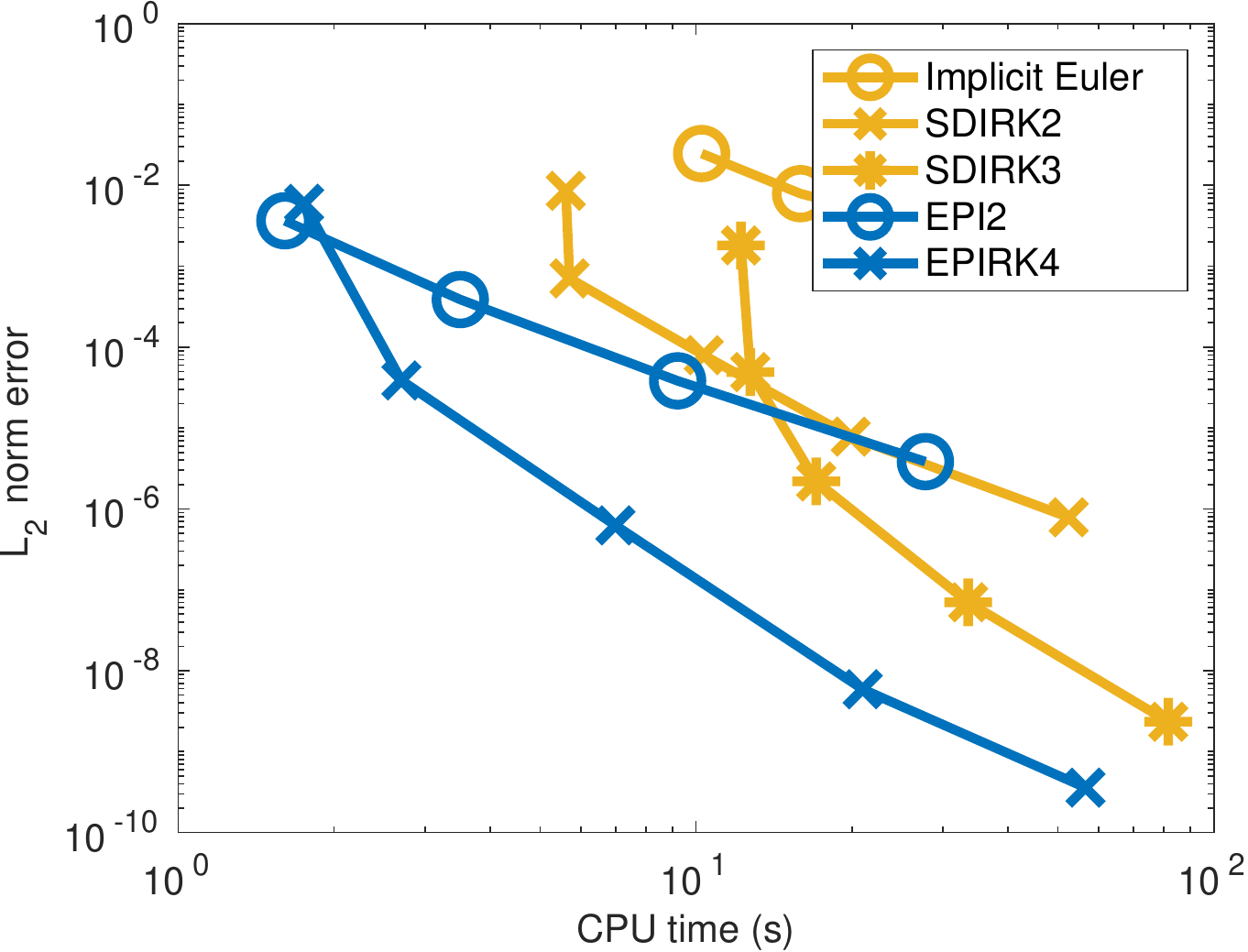}
    \caption{Implicit methods - Fine grid ($n_{elem} = 50^2$)}
  \end{subfigure}

  \caption{Precision diagrams for the 2d test problems with anisotropic diffusion ($\kappa = 10^{-2}, \alpha = 10^{-3}, \beta_1 = 0, \beta_2=10$).}
  \label{fig:prec_2d}
\end{figure}

\Cref{fig:prec_1d,fig:prec_2d} are showing the precision diagram (Error vs time to solution) for the 1d nonlinear diffusion and 2d anisotropic diffusion respectively. Each plot compare the performance of the exponential schemes to explicit (left) and implicit (right) methods on a coarse (top) and fine (bottom) grid. 
As expected, if we compare an explicit with an exponential scheme at the same order of accuracy, the time to solution for the explicit scheme will be lower as the cost per iteration is much cheaper with explicit scheme. 
However, as we can see from the previous section, explicit schemes quickly become unstable as the grid is refined and therefore can not be considered for stiff problems.
On the other hand, we can see that the performance of exponential schemes is better than the implicit methods on both test problems. Moreover, the gap between the implicit and exponential methods is increasing going from the coarse to the finner grid. This indicates that exponential methods scale better as the problem get stiffer.

\clearpage

\section{Conclusion and future work}
In this work, we introduced two nonlinear diffusion PDEs as a test problems for the RMHD equations. We used these problems to compare the performance of exponential integrators with traditional methods. These integrators and test problems were validated by looking at the order of convergence. Moreover, the results demonstrate that exponential integrators have the potential to be an efficient alternative to implicit methods for this kind of problems. 

In the future, we are planning on extending this work to more realistic test problems. We also want to investigate the addition of algebraic constraint as part of the solve. These constraint are of great interest for RMHD problems and need to be treated correctly and efficiently. 

\section{Acknowledgements}
Work for LLNL-TR-832157 was performed by LLNL under the auspices of the U.S. DOE under contract DE-AC52-07NA2734. The work was also supported by LLNL Laboratory Directed Research and Development project PLS-20-ERD-038 and by the NSF grants DMS-1720495 and DMS-2012875.
We would like to thank Chris Vogl and Milan Holec for their help, valuable discussions and work on the mhdex project. 

\pagebreak
\bibliographystyle{siam}
\bibliography{biblio.bib}

\begin{thebibliography}{1}

\bibitem{al-mohyComputingActionMatrix2011a}
{\sc A.~H. Al-Mohy and N.~J. Higham}, {\em Computing the {{Action}} of the
  {{Matrix Exponential}}, with an {{ Application}} to {{Exponential
  Integrators}}}, 33 (2011), pp.~488--511.

\bibitem{einkemmer2017a}
{\sc L.~Einkemmer, M.~Tokman, and J.~Loffeld}, {\em On the performance of
  exponential integrators for problems in magnetohydrodynamics}, 330 (2017),
  pp.~550--565.

\bibitem{gaudreault2022high}
{\sc S.~Gaudreault, M.~Charron, V.~Dallerit, and M.~Tokman}, {\em High-order
  numerical solutions to the shallow-water equations on the rotated
  cubed-sphere grid}, Journal of Computational Physics, 449 (2022), p.~110792.

\bibitem{gaudreaultKIOPSFastAdaptive2018}
{\sc S.~Gaudreault, G.~Rainwater, and M.~Tokman}, {\em {{KIOPS}}: {{A}} fast
  adaptive {{Krylov}} subspace solver for exponential integrators}, 372 (2018),
  pp.~236--255.

\bibitem{niesen2012a}
{\sc J.~Niesen and W.~Wright}, {\em Algorithm 919: A {{Krylov}} subspace
  algorithm for evaluating the $ \varphi$-functions appearing in exponential
  integrators}, 38 (2012), p.~22.

\bibitem{rainwater2016new}
{\sc G.~Rainwater and M.~Tokman}, {\em A new approach to constructing efficient
  stiffly accurate epirk methods}, Journal of Computational Physics, 323
  (2016), pp.~283--309.

\bibitem{tokmanEfficientIntegrationLarge2006}
{\sc M.~Tokman}, {\em Efficient integration of large stiff systems of {{ODEs}}
  with exponential propagation iterative ({{EPI}}) methods}, 213 (2006),
  pp.~748--776.

\bibitem{tokmanNewClassExponential2011}
\leavevmode\vrule height 2pt depth -1.6pt width 23pt, {\em A new class of
  exponential propagation iterative methods of {{Runge}} –{{Kutta}} type
  ({{EPIRK}})}, 230 (2011), pp.~8762--8778.

\bibitem{tokmanNewAdaptiveExponential2012}
{\sc M.~Tokman, J.~Loffeld, and P.~Tranquilli}, {\em New {{Adaptive Exponential
  Propagation Iterative Methods}} of {{ Runge--Kutta Type}}}, 34 (2012),
  pp.~A2650--A2669.

\end{thebibliography}
 
\end{document}